# Application of the equal dissipation rate principle to automatic generation of strut-and-tie models

B. Novák & S. Ananiev
*Institute of Lightweight Structures and Conceptual Design*
*University of Stuttgart, Pfaffenwaldring 7, 70569 Stuttgart, Germany*
*E-mail: balthasar.novak@ilek.uni-stuttgart.de, sergey.ananiev@ilek.uni-stuttgart.de*

ABSTRACT: This work presents an extended formulation of maximal stiffness design, within the framework of the topology optimization. The mathematical formulation of the optimization problem is based on the postulated principle of equal dissipation rate during inelastic deformation. This principle leads to the enforcement of stress constraints in domains where inelastic deformation would occur. During the transition from the continuous structure to the truss-like one (strut-and-tie model) the dissipation rate is kept constant using the projected gradient method. The equal dissipation rate in the resulting truss and in the original continuous structure can be regarded as an equivalence statement and suggests an alternative understanding of physical motivation behind the strut-and-tie modeling. The performance of the proposed formulation is demonstrated with the help of two examples.

## 1 MOTIVATION

### 1.1 Why strut-and-tie models?

Despite the rapid development of numerical methods in structural mechanics, the fundamentals of conceptual design did not change much. A good engineer can be distinguished by his or her ability to model the load-bearing mechanism of a structure using simplified models. This ability can be gained with years of practical experience and has a decisive meaning in prevention of crude errors in the design process. Such errors can be hardly considered during numerical simulation.

This is especially true for the design of concrete structures. The material properties of concrete have a very complicated nature and are still an active research field. There are many different theories like plasticity, continuum damage, microplane models, fracture mechanics or discrete lattice models that can be applied to the modeling of concrete. Nevertheless concrete has been successfully used for more than one hundred years in civil engineering. Obviously, this was only possible due to the simple and robust approaches developed for its modeling. One of these approaches is to replace the continuum structure with a truss-like one, where stress state is nearly one-dimensional. Pioneering contributions in this field were made by Ritter (1899) and Mörsch (1902). In these works, truss models were developed for simple beam structures. The idea has proved to be powerful and became quickly a common practice.

A solid attempt to generalize this approach to discontinuity regions, where the Bernoulli hypothesis is not valid, was undertaken by Marti (1985) and Schlaich et al. (1987). It was proposed to develop strut-and-tie models, which can be even kinematical, but which closely follow the elastic stress distribution. If one achieves that, then a realistic prediction of load-bearing capacity of the whole structure can be made.

Nowadays, after years of research and practical application (Schlaich & Schäfer 2001), the strut-and-tie modeling offers a valuable design tool, which is recommended for use by several national Design Codes (Reineck 2002).

### 1.2 Why automatic generation?

A direct understanding of optimization procedure as a tool for automatic generation of strut-and-tie models (STM) would contradict the considerations presented above. Indeed, it was stressed that the main value of STMs is the clarity and simplicity of the approach, where the qualification and the intuition of designer plays the central role. The "automatization tool" can potentially devaluate this role to the responsibility of correct input of boundary conditions. This is definitely not the purpose of this work.

On the other hand, the uniqueness of the models is still an open question (Reineck 2002). It is not clear according to which criterion two different models for the same discontinuity regions can be compared. This uncertainty constitutes a serious ob-stacle for the broad acceptance of STMs in the engineering praxis. The development of such criterion is the main objective of this work.

## 2 EXISTING APPROACHES

### 2.1 *Integration of stress fields*

The engineering literature (Schlaich & Schäfer 2001, Reineck 2002) recommends the use of elastic stress distributions as basis for STM design. The struts or ties appear here as the center line of the integrated stress fields. The best model is those, which fulfills the following energy criteria:

$$\min \sum_{i=1}^{N} F_i L_i \varepsilon_i \qquad (1)$$

where $F_i$ is the force in the truss member i, $L_i$ is its length, $\varepsilon_i$ is its averaged strain and N is the number of truss members.

Figure 1 shows the classical application of this approach to the concrete deep beam of Leonhardt & Walther (1966).

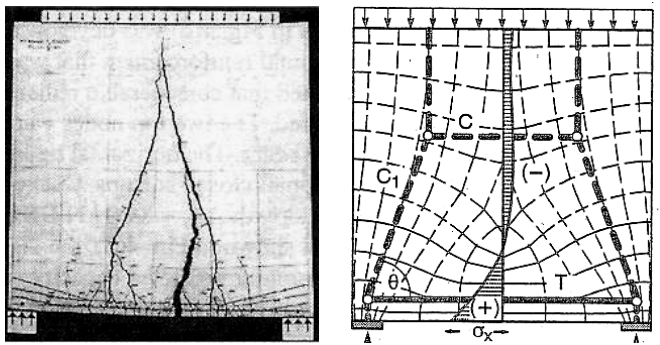

Figure 1. The deep beam of Leonhardt & Walther (1966) and the corresponding STM from Schlaich & Schäfer (2001).

There are several restrictions of this approach, which can be summarized as follows:
– the usage of elastic stress distribution is not easy to justify. In discontinuity regions the stress state is usually highly inelastic, with partial dissipation of energy;
– it is also difficult to justify the replacement of continuous structure by some truss-like model. The application of the lower bound theorem of the theory of plasticity is not entirely correct because of the assumed unlimited material ductility. Furthermore the inelastic deformation in concrete due to cracks and internal damage cannot be understood as plastic deformations;
– special attention has to be paid to the boundary conditions. It is not completely correct to replace the homogeneous pressure on the bottom of the beam by concentrated loads (Fig. 1). This introduces a new discontinuity, which is not present in the original structure.

These inconsistencies often lead to unrealistic STM designs. For example, the increasing or decreasing of the loading does not affect the geometry of the models. It is physically more adequate to assume that the STMs must be different for each stage of loading. This point is often argued in the engineering literature (Reineck 2002).

### 2.2 *Maximal damage resistance design*

The STM is not the only way to design the reinforcement in concrete elements. In the more general approach the optimal reinforcement is understood as those, which minimizes the dissipated energy during the whole loading history. Several variations of this problem were considered by Ramm et al. (1998), Hammer (2000), Barthold & Firuziaan (2000) and Ananiev (in prep.). Within this framework, the structure is treated as a continuum, where a suitable inelastic material model is used for dissipative processes. The amount of reinforcement in each element is taken as a design variable. The obvious advantage of this approach is its general character: the structure can be optimized for an arbitrary loading history. Consistent linearization of the material model and the history dependent calculation of gradients constitute the most part of the work. (Michaleris et al. 1994)

Its disadvantage has a practical nature. It is not clear how the "clouds" of optimal reinforcement can be transformed to strut and ties. The equilibrated stress state achieved at the end of the loading history is unique for the continuum. Any attempt to concentrate reinforcement and concrete in a truss-like structure will destroy it. That means that the new equilibrium state may have nothing to do with the original structure.

Figure 2 shows the typical optimal reinforcement and load-displacement diagram obtained by this approach. Due to the restriction on available amount of reinforcement, the load-bearing capacity of optimal structure is limited to $F^{MAX}$.

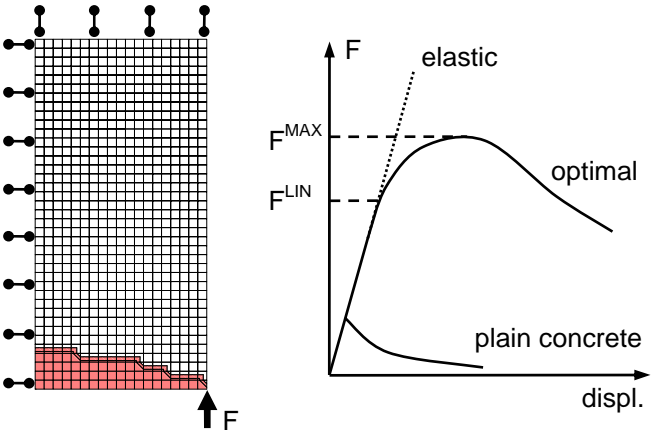

Figure 2. Optimal reinforcement and typical load-displacement diagram for maximal damage resistance design (from Ananiev, in prep.).

As can be seen from Figure 2, this approach can not be directly used for automatic generation of strut-and-tie models. However, this approach delivers important information about the areas, where the optimal reinforcement must be placed to minimize the dissipation of energy. This result will be used implicitly in our work.

## 2.3 Optimization of discrete truss-like structures

The truss with all possible connections of nodes can be taken directly as the base structure. This type of optimization problem was considered by Ramm et al. (1998) and Ali & White (2001), where the cross section's area of truss member was taken as a design variable.

This approach has a theoretical limitation. One has to remember that the main objective of STM design is an optimal reinforcement in *continuous* concrete structures. It is not trivial to prove that the stress state in the discrete truss will converge to the continuous stress state for arbitrary material parameters. This task will become even more difficult if the inelastic material behavior is taken into account. Without such a proof, it is not possible to assure that the optimized truss has something to do with the optimal reinforcement of the original continuous concrete structure.

## 3 PRINCIPLE OF EQUAL DISSIPATION RATE

### 3.1 Basic assumption

It has become clear from the previous discussion that this work has two main objectives, namely
- to develop a method for optimal reinforcement design in continuous concrete structures;
- to represent the results in a truss-like form, where the stress state is one-dimensional.

At the first glance, these objectives are self-contradictory. Indeed, only under very special assumptions, it seems to be possible to achieve them.

In order to justify these assumptions, one has to return to the objective of maximal damage resistance design, which aims a complete compensation of the dissipated energy in concrete. In other words, it attempts to follow the elastic deformation path as long as possible. Therefore, for the case of a moderate loading, like $F^{LIN}$ in Figure 2, an *elastic behavior* can be assumed. This assumption allows the use of the reactions obtained from linear-elastic calculations, leading to a simplified treatment of statically indeterminate systems.

Another important consequence of this assumption becomes obvious in the case of proportional loading. Optimal reinforcement can be generated in just one iteration: all concrete elements with inelastic stress state must be replaced with steel elements of equal stiffness.

### 3.2 Formulation of the principle

To transform the "clouds" of reinforcement into a truss-like structure we have to define the criterion according to which two structures with different stress states can be understood as equivalent.

To this end, we state the following equivalence principle: two structures are equivalent if they have the same dissipation rate.

Its mathematical formulation reads:

$$\mathcal{D}_1 = \mathcal{D}_2$$
$$\mathcal{D}_i = \int_{V_i} \boldsymbol{\sigma} : \dot{\boldsymbol{\varepsilon}}^{in} dv \quad \text{with} \quad \dot{\boldsymbol{\varepsilon}}^{in} = \lambda \frac{\partial f}{\partial \boldsymbol{\sigma}}, \qquad (2)$$

where $\boldsymbol{\sigma}$ is the elastic stress tensor, $\dot{\boldsymbol{\varepsilon}}^{In}$ is the inelastic rate of strain tensor, V is the volume of the structural or finite element, f is the loading surface and $\lambda$ is a Lagrange multiplier.

The inelastic strains are understood here in the general sense proposed by Meschke et al. (1998): if the dissipative process is plastic, than they are real (plastic) strains, if the inelasticity has different nature (e.g. damage), then they are fictitious and the real quantities which grow are damage variables (see also Ananiev & Ožbolt 2004). The rates of inelastic strains are determined using the associated flow rule.

The form of the loading surface, which is required by the flow rule in Equation 2, constitutes the second important assumption in this work. Taking into account that our primary goal is to obtain the optimal reinforcement in the damaged concrete areas, we impose the Rankine tension cut-off criterion ($\sigma_I \leq \tau_0$). This leads to the following dissipation rate in a structure:

$$\mathcal{D}_i = \int_{V_i} \lambda \sigma_I dv, \qquad (3)$$

where $\sigma_I$ is the first principal stress.

For simplicity reasons it is also assumed that the compressive strength of concrete is unlimited.

## 3.3 *Stresses as inequality constraint*

The introduction of stress constraint is not a new idea. Duysinx & Bendsøe (1998) and Ramm et al. (1998) have already proposed to enrich the optimization formulation with *inequality* constraint in the following form:

$$f(\boldsymbol{\sigma}) \leq (\rho_e)^p \tau_0, \quad (4)$$

where $f(\boldsymbol{\sigma})$ is the loading function (e.g. von Mises), $\tau_0$ is the equivalent yield stress and p is a penalty parameter.

There are many engineering problems, where this formulation is meaningful. However, this is not the case for reinforcement design. Our objective is not to avoid the violation of stress constraints – they are already violated – but to compensate the loss of energy in the structure.

The physical inconsistency of this formulation can be illustrated by the simple example with two rods, which is shown in the Figure 3. The objective here is the structure with maximal stiffness and the cross section's areas of rods are taken as design variables. The amount of material remains constant: 3AL.

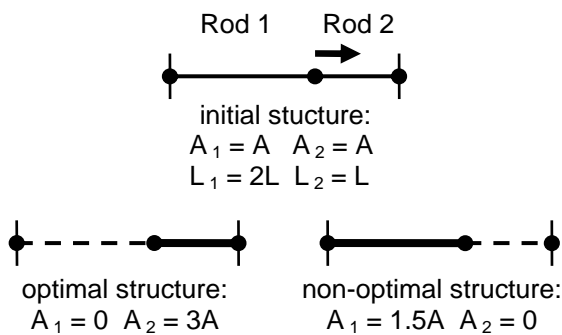

Figure 3. Example of an optimal structure, which does not violate stress constraint, but cannot be accepted as an optimal reinforcement design.

Obviously, the shortest rod represents an optimal structure, because this is the stiffest model. According to it, there is no need of reinforcement, which is not correct from the engineering point of view, because of the uncertainty of the concrete tensile strength. Such an "optimal design" could lead to the complete cracking of the left rod. In the sense of our equivalence principle (Eq. 2) these two structures (initial and optimal) are not equivalent. The equivalence is kept if the stress in the first rod remains *constant* during optimization. For this example this means that the initial structure will not be changed at all (the number of constraints is equal to the number of design variables).

## 4 FORMULATION OF THE OPTIMIZATION PROBLEM

### 4.1 *Definition of design variables*

Due to the theoretical limitation imposed by discrete truss-like structures (Sec.2.3) the discretized continuum is chosen as the base structure. The density of material in each element serves as a design variable. This type of optimization problem belongs to the continuum topology optimization. An excellent state-of-the-art review of this field can be found in the recent book of Bendsøe & Sigmund (2003).

The topology optimization was already applied to the automatic generation of STM by Ramm et al. (1998) and Liang et al. (2001). The present work can be regarded as an extension of their ideas.

### 4.2 *Objective function and the mass constraint*

To assure that the optimization problem leads to the concentration of uniformly distributed material in truss-like structures, the following two components of the formulation are essential: the objective function, which is the stored elastic energy and the mass constraint.

$$\min \int_V \frac{1}{2} \boldsymbol{\sigma} : \boldsymbol{\varepsilon} \, dv \quad \text{s.t.} \quad \int_V \rho \, dv = M, \quad (5)$$

where $\rho$ is the material density and M is the total amount of material to be redistributed.

The minimization of the stored elastic energy means that the optimal structure is the stiffest one. It is interesting to note that this optimality criterion have been used intuitively in the engineering literature for a long time. For example, for one dimensional stress state we have (dv = da·dl):

$$\int_V \boldsymbol{\sigma} : \boldsymbol{\varepsilon} \, dv = \int_V \boldsymbol{\sigma} da : \boldsymbol{\varepsilon} dl = \sum_{i=1}^{N} F_i L_i \varepsilon_i \quad (6)$$

So, the objective function complies with the quality criterion for STM shown in Equation 1. Therefore, the structures generated using this objective function will automatically fulfil the required criteria of Schlaich & Schäfer (2001).

Imposing of mass constraint leads to the optimal structure, which has only simple one dimensional stress states. The material will concentrate in truss-like structures, because the existence of any bending in structure would be ineffective. The mathematical explanation of this fact will be given in the next section while examining the structure of projected gradient method.

### 4.3 *Solution method*

In the works on topology optimization by Bendsøe & Sigmund (2003), the heuristic optimality criteria method is often used for the numerical solution. It was shown recently (Ananiev 2005) that with definition of Lagrange multipliers proposed by Hestenes (1975), this method is equivalent to the projected gradient method. The central idea was to introduce the multipliers not at the optimum, but during projection of the objective function's gradient onto tangential space of active constraints. Aside from its rigorous mathematical structure, this method is relatively simple to implement and to analyze. Already in the first iteration, this method can show if the chosen optimization formulation will lead to the desired optimum.

The projected gradient is constructed as a linear combination of objective function's gradient and gradients of active constraints:

$$\vec{d} = -\nabla f + \sum_{k=1}^{NH} \lambda_k \nabla h^k, \quad (7)$$

where $\nabla f$ is the gradient of objective function f, $\nabla h^k$ is the gradient of the active constraint $h^k$, $\lambda_k$ is the corresponding Lagrange (Hestenes) multiplier and NH is the number of active constraints.

The multipliers $\lambda_k$ are found from the orthogonality condition between vector $\vec{d}$ and the gradients of all active constraints ($\nabla h^k$, k = 1..NH). This condition means that an infinitesimal small change of design variables along projected gradients will not violate the active constraints.

$$\begin{cases} \vec{d} \cdot \nabla h^I = 0 \\ \quad \ldots \\ \vec{d} \cdot \nabla h^{NH} = 0 \end{cases} \Rightarrow \vec{\lambda} = \left[ \mathbf{H} \mathbf{H}^T \right]^{-1} \cdot \left[ \nabla f \cdot \mathbf{H} \right] \quad (8)$$

The exam of the projected gradient's structure allows explaining why the mass constraint was introduced in Equation 5. With its help, the projected gradient has always positive and negative components. For example, in the case of structured grid with elements of equal volume, the Lagrange (Hestenes) multiplier $\lambda$ is equal to the mean value of the objective function's gradient.

$$\lambda = \frac{\sum_{e=1}^{N} f_{,\rho_e} \cdot h^{II}_{,\rho_e}}{\sum_{e=1}^{N} h^{II}_{,\rho_e} \cdot h^{II}_{,\rho_e}} = \frac{\sum_{e=1}^{N} f_{,\rho_e} \cdot 1}{N} \quad (9)$$

This means that the density in the most stressed elements will grow, while in the others it will decrease. Thus, the concentration of material in chains is guaranteed. Figure 4 illustrates the change in structure of the objective function's gradient after the projection.

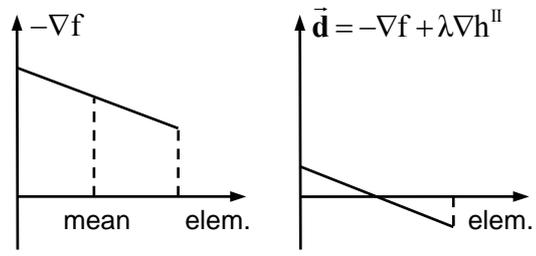

Figure 4. The structure of projected gradient if only mass constraint is active.

After the projected gradient is built, an update of design variables is done using the simple formula:

$$\vec{\rho}^{,n+1} = \vec{\rho}^{,n} + \gamma \vec{d}^{,n}, \quad (10)$$

where $\vec{\rho}^{,n+1}$ is the new value of the densities, $\vec{\rho}^{,n}$ is the current value of the densities, $\vec{d}^{,n}$ is the projected gradient, n is the iteration number and $\gamma$ is the size of the update vector (usually, $\gamma \ll 1$).

## 5 CALCULATION OF THE GRADIENTS

### 5.1 *Discretized form of the optimization problem*

After discretization with finite elements, the complete formulation of the extended optimization problem discussed in Section 4 looks as follows:

$$\min_{\vec{\rho}} \quad f: \vec{u}^T \cdot \mathbf{K}(\vec{\rho}) \cdot \vec{u}$$

subject to (11)

$$\begin{cases} h^I \quad : \mathbf{K}(\vec{\rho}) \cdot \vec{u} = \vec{p} \\ h^{II} \quad : \vec{\rho} \cdot \vec{v} = \sum_{e=1}^{N} \rho_e v_e = M \\ h^{IIIe} : 0 \leq \rho_e \leq 1, \; e = 1..N \\ h^{IVd} : \sum_{e=1}^{NDE} \left\{ \sum_{i=1}^{4} \sigma_{Ie}(x_i, y_i) \det \mathbf{J}_e \right\} = const, \; d = 1..ND \end{cases}$$

where f is the objective function, $h^i$ (i = I..IV) are the constraints, $\vec{u}$ is the displacement vector of the whole structure, $\mathbf{K}(\vec{\rho})$ is the stiffness matrix of the whole structure, $\vec{\rho}$ is the vector of the element densities, $\vec{p}$ is the vector of the external loads, $\vec{v}$ is the vector of the element volume, M is the mass of the material, $\rho_e$ is the material density in element e, $v_e$ is the volume of element e, d =1..ND is the number of inelastic domains, NDE is the number of elements in the each inelastic domain, $\sigma_{Ie}(x_i, y_i)$ is the first principal tensile stress in the Gauss point i, where i = 1..4 and $\det \mathbf{J}_e$ is the Jacobian.

Formally, the optimization problem stated in the system of equations (11) depends on two variables, namely the densities ($\vec{\rho}$) and the displacements ($\vec{u}$). Due to the equilibrium constraint ($h^I$), the primal design variables are the element material densities and the displacements depend on them in implicit way. Therefore, all gradients will be built with respect to the densities.

## 5.2 The gradient of the objective function

The application of the chain rule to the objective function in Equation 11 furnishes

$$\frac{\partial f}{\partial \rho_e} = 2\frac{\partial \vec{u}^T}{\partial \rho_e} \cdot \mathbf{K}(\vec{\rho}) + \vec{u}^T \cdot \frac{\partial \mathbf{K}(\vec{\rho})}{\partial \rho_e} \cdot \vec{u}. \quad (12)$$

The sensitivity of the displacements are obtained from the equilibrium constraint, as follows

$$\frac{\partial}{\partial \rho_e}\{\mathbf{K}(\vec{\rho})\cdot\vec{u}-\vec{p}\}=0 \rightarrow \frac{\partial \vec{u}}{\partial \rho_e} = -\mathbf{K}^{-1}\cdot\left[\frac{\partial \mathbf{K}}{\partial \rho_e}\cdot\vec{u}\right] \quad (13)$$

Substituting this result into Equation 12 and taking into account the local character of the dependence of the stiffness matrix on the density in an element, it leads to

$$\frac{\partial f}{\partial \rho_e} = -\vec{u}^T \cdot \frac{\partial \mathbf{K}(\vec{\rho})}{\partial \rho_e} \cdot \vec{u} = -\vec{u}_e^T \cdot \frac{d\mathbf{k}_e}{d\rho_e} \cdot \vec{u}_e$$

$$\frac{d\mathbf{k}_e}{d\rho_e} = \frac{p}{\rho_e}\mathbf{k}_e^0 \quad (14)$$

where $\vec{u}_e$ is the displacement vector of element e, $\mathbf{k}_e^0$ is the initial stiffness matrix of element and p is the penalty parameter.

The dependence of the element stiffness matrix on the actual value of the density is described using the Simple Isotropic Model with Penalization (SIMP), $\mathbf{k}_e = (\rho_e)^p \mathbf{k}_e^0$ (Bendsøe & Sigmund 2003).

## 5.3 The gradient of the constant mass constraint

Due to the linearity of the mass constraint, its gradient is a constant vector equal to the vector of element volumes:

$$\frac{\partial h^{II}}{\partial \vec{\rho}} = \vec{v} \quad (15)$$

## 5.4 The gradient of the boundary constraint

Because of the inequality character of this constraint, one has to check if it is active or not. It is done using an iterative procedure, where the densities are updated under the assumption that there are no active boundary constraints. If after this trial update some of them are violated, than the corresponding densities are fixed at the old values and the new projected gradient is built. These trial updates are repeated until the number of active constraints does not change (Active Set Method, see e.g. Luenberger 1989).

The gradient of the active boundary constraint is a Kronecker-delta vector:

$$\frac{\partial h^{IIIs}}{\partial \rho_e} = \delta_{se} = \begin{cases} 1, & \text{if}(s=e) \\ 0, & \text{if}(s \neq e) \end{cases} \quad (16)$$

## 5.5 The gradient of the dissipation rate

Following the standard scheme of the Finite Element Method, the dissipation rate from Equation 3 is evaluated numerically at element level using Gauss quadrature (Eq. 11). All rates in elements belonging to the same inelastic domain are added together. It is important to understand that the dissipation rate is held constant not in each inelastic element, but in some domain of such elements. This allows the "concentration" of reinforcement in ties. It is also important to mention that during the optimization, elements cannot change their type even if their stress state becomes elastic. For example, in the Leonhardt & Walther deep beam (Fig. 2) there is only one inelastic domain.

The application of the chain rule requires the calculation of the sensitivity of three quantities: the Lagrange multipliers, the first principal stresses and the Jacobian.

According to our assumption of elastic behavior at the structural level (Sec.3.1) there is no actual inelastic deformation. One can regard such stress state as an elastic predictor, which is known from the computational inelasticity (Simo & Hughes 1998). The introduced dissipation constraint concerns the speed of energy loss if we *hypothetically* allow the inelasticity in concrete. Due to this reason, the sensitivity of the Lagrange multiplier is a zero vector (all $\lambda = 1$).

The change of material density does not affect the form of finite element. Therefore, the sensitivity of Jacobian is a zero vector as well.

$$\frac{\partial \lambda_e}{\partial \vec{\rho}} = \mathbf{0}, \quad \frac{\partial \det \mathbf{J}_e}{\partial \vec{\rho}} = \mathbf{0}, \quad \frac{\partial \sigma_{Ie}}{\partial \vec{\rho}} = ? \quad (17)$$

Before we proceed with sensitivity of the stresses, it is necessary to clarify if the tensile strength of concrete and/or stresses have to be multiplied with penalized density: $(\rho_e)^p$ similar to the Young modulus in SIMP model (Sec. 3.3, 5.2).

One of the further advantages of our (linear) elasticity assumption is the fact that stresses in the optimized structure do not depend on the initial value of material density (all $\rho_e^0 = \rho^{Start}$). To prove this, let us first consider the change in the global displacement vector:

$$\mathbf{K}(\vec{\rho}^{\text{Start}}) \cdot \vec{u} = \vec{p} \rightarrow (\rho^{\text{Start}})^p \mathbf{K}^{,0} \cdot \vec{u} = \vec{p} \rightarrow$$

$$\vec{u} = \frac{1}{(\rho^{\text{Start}})^p} [\mathbf{K}^{,0}]^{-1} \cdot \vec{p} = \frac{1}{(\rho^{\text{Start}})^p} \vec{u}^{,0}, \quad (18)$$

where $\mathbf{K}^{,0}$ is the initial stiffness matrix of the whole structure, $\vec{u}^{,0}$ is the corresponding initial displacement vector, $\vec{u}$ is the current displacement vector.

Second, we insert this scaled displacement vector into the expression for the calculation of the stress at a Gauss point.

$$\vec{\sigma}_e = (\rho^{\text{Start}})^p \mathbf{D}^{,0} \mathbf{B} \cdot \vec{u}_e = (\rho^{\text{Start}})^p \mathbf{D}^{,0} \mathbf{B} \cdot \frac{\vec{u}_e^{,0}}{(\rho^{\text{Start}})^p} \quad (19)$$

$$\vec{\sigma}_e = \vec{\sigma}_e^{,0}$$

where $\vec{\sigma}_e$ is the vector representation of the stress tensor, $\mathbf{D}^{,0}$ is the matrix representation of initial elasticity tensor, $\mathbf{B}$ is the strain-displacement matrix and $\vec{u}_e$ is the vector of element displacements.

The initial value of density is not present in the Equation 19 and therefore the stresses preserve their initial values. This gives a physical motivation to perform the sensitivity analysis of the *actual* values of principal stresses during the current optimization iteration. In other words, the dissipation rate is kept constant during the whole optimization history.

Now we are in position to build the gradients of the first principal stresses. For two-dimensional stress state we have:

$$\frac{\partial \sigma_I}{\partial \rho_k} = \frac{\partial \sigma_I}{\partial \sigma_{xx}} \frac{\partial \sigma_{xx}}{\partial \rho_k} + \frac{\partial \sigma_I}{\partial \sigma_{yy}} \frac{\partial \sigma_{yy}}{\partial \rho_k} + \frac{\partial \sigma_I}{\partial \sigma_{xy}} \frac{\partial \sigma_{xy}}{\partial \rho_k} \quad (20)$$

with

$$\sigma_I = \frac{1}{2}(\sigma_{xx} + \sigma_{yy}) + \frac{1}{2}\sqrt{(\sigma_{xx} + \sigma_{yy})^2 + 4(\sigma_{xx}\sigma_{yy} - \sigma_{xy}^2)} \quad (21)$$

The sensitivities of stress tensor components in Equation 20 are determined using the linear elasticity relations (SIMP model):

$$\begin{bmatrix} \sigma_{xx,\rho_k} \\ \sigma_{yy,\rho_k} \\ \sigma_{xy,\rho_k} \end{bmatrix} = \begin{cases} (\rho_e)^p \mathbf{D}^{,0} \mathbf{B} \cdot \frac{\partial \vec{u}_e}{\partial \rho_k}, & \text{if } (k \neq e) \\ p(\rho_e)^{p-1} \mathbf{D}^{,0} \mathbf{B} \cdot \vec{u}_e + (\rho_e)^p \mathbf{D}^{,0} \mathbf{B} \cdot \frac{\partial \vec{u}_e}{\partial \rho_k}, & \text{if } (k = e) \end{cases} \quad (22)$$

The sensitivities of element's displacements are already known at the structural level (Eq. 13).

Finally, we obtain the gradient of the dissipation rate in some inelastic domain:

$$\frac{\partial h^{IVd}}{\partial \rho_k} = \sum_{e=1}^{NDE} \left\{ \sum_{i=1}^{4} \frac{\partial \sigma_{Ie}(x_i, y_i)}{\partial \rho_k} \det \mathbf{J}_e \right\} \quad (23)$$

# 6 EXAMPLES

## 6.1 *Deep beam of Leonhardt & Walther*

As already mentioned in Section 2.1, special attention has to be paid to the boundary conditions in this example.

The discontinuity in this case is the support of the deep beam. It is not correct to replace the homogenous pressure at the top of the deep beam with two statically equivalent concentrated loads, because this would introduce a new discontinuity. We make use of the fact that in the Bernoulli region not only the pressure is homogenous, but also the displacements. This can be modeled by distributed supports at the top of the deep beam. The statically equivalent force is applied at the bottom (Fig. 5). The additional advantage of this approach is that the distance between compression struts will be determined automatically. As it can be seen from Figures 5-7, it is always larger than half of the deep beam's width (compare with Figure 1).

Three different optimal strut-and-tie models were generated using the proposed extended formulation (Eq. 11). In all cases, the same statically equivalent force is applied, but the tensile strength of concrete varies. To simplify the understanding of the results, each figure contains also the inelastic domains, where the dissipation rate is kept constant.

The results are shown in the Figures 5-7 followed by the diagram (Fig. 8) showing the evolution of the objective function: stored elastic energy. As it could be expected, the stiffest of these three examples is simple vertical compression strut (Fig. 7), which is the optimal strut-and-tie model, if all concrete elements are elastic.

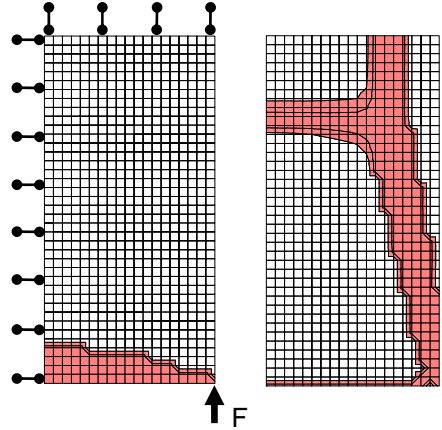

Figure 5. Inelastic domain and the corresponding optimal strut-and-tie model using tensile strength = 1f.

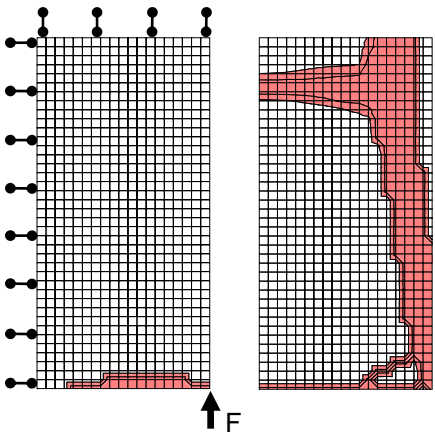

Figure 6. Inelastic domain and the corresponding optimal strut-and-tie model using tensile strength = 2f.

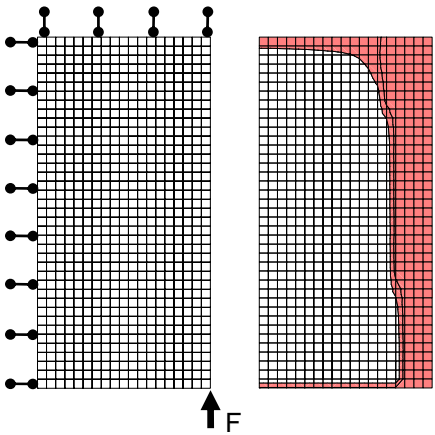

Figure 7. Inelastic domain and the corresponding optimal strut-and-tie model using tensile strength = 10f (all elements are elastic).

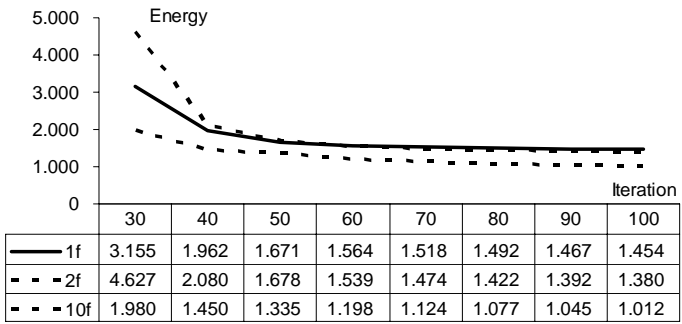

| | 30 | 40 | 50 | 60 | 70 | 80 | 90 | 100 |
|---|---|---|---|---|---|---|---|---|
| 1f | 3.155 | 1.962 | 1.671 | 1.564 | 1.518 | 1.492 | 1.467 | 1.454 |
| 2f | 4.627 | 2.080 | 1.678 | 1.539 | 1.474 | 1.422 | 1.392 | 1.380 |
| 10f | 1.980 | 1.450 | 1.335 | 1.198 | 1.124 | 1.077 | 1.045 | 1.012 |

Figure 8. Convergence history of stored elastic energy in the Leonhardt & Walther deep beam.

## 6.2 *Deep beam with opening*

In contrast to the previous example, this one is more of academic interest. It allows demonstrating the performance of the proposed formulation in case of several inelastic domains. Due to the double symmetry of the problem, only one fourth of the deep beam is considered.

Figure 10 shows that for the chosen loading conditions there are three inelastic domains. In the sense of our extended formulation (Eq. 11) there are three constraints of type $h^{IV}$. The number of elements in each domain remains constant during the optimization process.

For illustration purposes, the same deep beam was also optimized with no inelastic domains (Fig.12). Similar to the Leonhardt & Walther example, this optimum represents the stiffest STM for the given geometry (Fig.13).

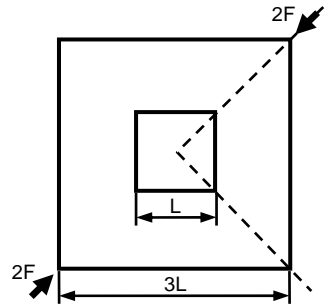

Figure 9. Concrete deep beam with opening.

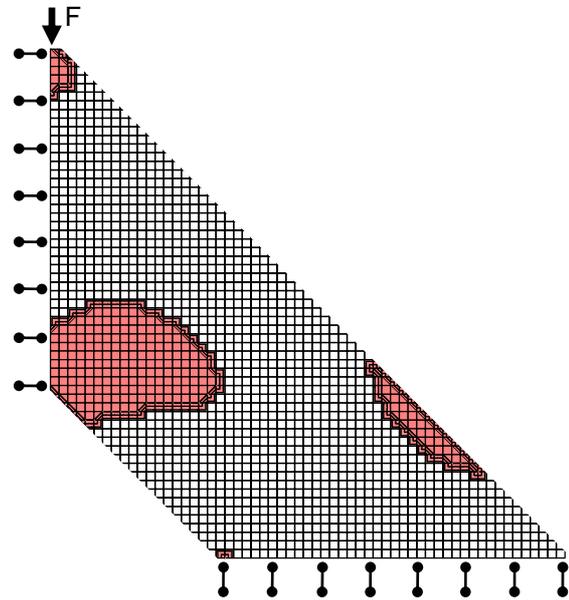

Figure 10. One fourth of the deep beam with three inelastic domains (tensile strength = 1f).

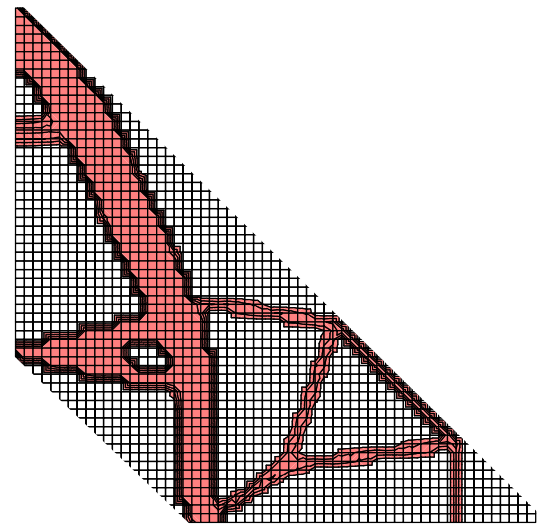

Figure 11. Optimal strut-and-tie model for three inelastic domains.

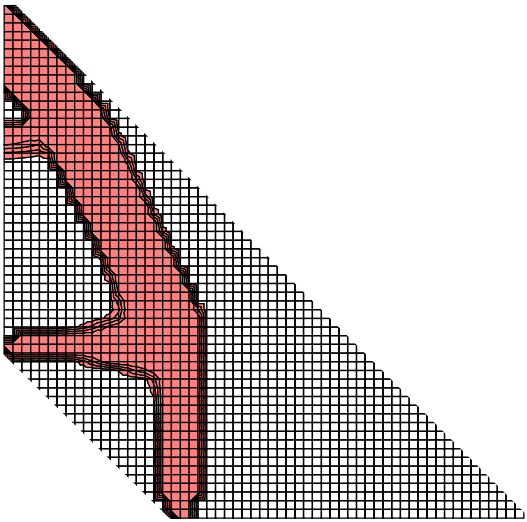

Figure 12. Optimal strut-and-tie model using tensile strength = 10f (all elements are elastic).

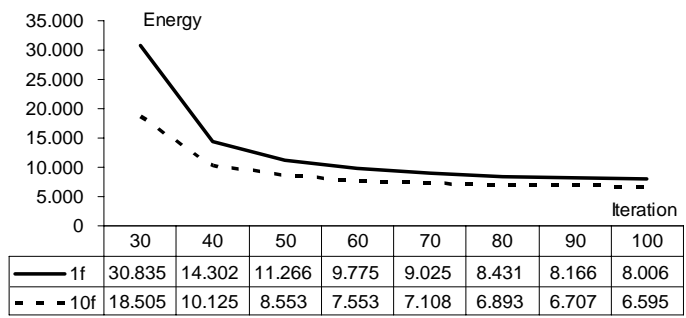

Figure 13. Convergence history of stored elastic energy in the deep beam with opening.

## 7 CLOSURE

The proposed extension of the maximal stiffness design is based on the assumption of complete compensation of dissipated energy at the structural level. This allows an effective formulation of the optimization problem and leads to reasonable strut-and-tie models. However, for practical purposes this assumption is obviously too strong, because of the inelastic behavior of reinforced concrete elements. The generalization of the proposed formulation also for the inelastic range represents a challenging problem for future research.

## ACKNOWLEDGEMENTS

A financial support of this work from the German Research Foundation (www.dfg.de) is gratefully acknowledged.

## REFERENCES

Ali, M.A. & White, R.N. 2001. Automatic generation of truss model for optimal design of reinforced concrete structures. *ACI Structural Journal* 98(4): 431-442.

Ananiev, S. 2005. On equivalence between optimality criteria and projected gradient methods with application to topology optimization problem. *Multibody System Dynamics* 13: 25-38.

Ananiev S. (in prep.). Zur nichtlinearen Optimierung in quasi-spröden Materialen am Beispiel von Stahlbeton. Ph.D. Thesis, University of Stuttgart.

Ananiev, S. & Ožbolt, J. 2004. Plastic-damage model for concrete in principal directions. In: Li V. et al. (eds.). *FraMCoS-5 Fracture Mechanics of Concrete and Concrete Structures*. Colorado, USA.

Barthold, F.-J. & Firuziaan, M. 2000. Optimization of hyperelastic materials with isotropic damage. *Structural and Multidisciplinary Optimization* 20: 12–21.

Bendsøe, M.P. & Sigmund, O. 2003. *Topology Optimization. Theory, Methods and Applications*. Berlin: Springer.

Duysinx, P. & Bendsøe, M.P. 1998. Topology optimization of continuum structures with local stress constrains. *International journal for numerical methods in engineering* 43: 1453-1478.

Hammer, V.B. 2000. Optimization of fibrous laminates undergoing progressive damage. *International journal for numerical methods in engineering* 48: 1265-1284.

Hestenes, M.R. 1975. *Optimization theory. The finite dimensional case*. New York: John Wiley & Sons.

Leonhardt, F. & Walther, R. 1966. Wandartiger Träger. *Deutscher Ausschuss für Stahlbeton* 178. Berlin: W. Ernst & Sohn.

Liang, Q.Q., Xie, Y.M., Steven, G.P. 2001. Generating optimal strut-and-tie models in prestressed concrete beams by performance-based optimization. *ACI Structural Journal* 98(2): 226-232.

Luenberger D. (1989). *Linear and nonlinear programming*. Massachusetts: Addison-Wesley Publishing Inc..

Marti, P. 1985. Basic tools of reinforced concrete beam design. *ACI Structural Journal* 82(1): 46-56.

Meschke, G., Lackner, R., Mang, H. 1998. An anisotropic elastoplastic-damage model for plain concrete. *International Journal of Numerical Methods in Engineering* 42: 703-727.

Michaleris, P., Tortorelli, D.A., Vidal, C.A. 1994. Tangent operators and design sensitivity formulations for transient nonlinear coupled problems with applications to elastoplasticity. *International Journal for Numerical Methods in Engineering*, 37(14): 2471-2501.

Mörsch, E. 1902. *Der Eisenbetonbau, seine Theorie und Anwendung*. Stuttgart: K. Wittwer. (English translation: Concrete steel construction. New York: McGraw-Hill, 1909).

Ramm, E., Maute, K., Schwarz, S. 1998. Conceptual design by structural optimization. In: de Borst R., Bićanić N., Mang H., Meschke G. (eds.). *Proceedings of EURO-C 1998 conference on computational modeling of concrete structures*: 879-896. Rotterdam: A.A. Balkena.

Reineck, K.-H. (eds) 2002. Examples for the design of structural concrete with strut-and-tie models. *Publication of ACI Committee 445: shear and torsion*. Farmington Hills: American Concrete Institute.

Ritter, W. 1899. Die Bauweise Hennebique. *Schweizerische Bauzeitung* XXXIII(7).

Schlaich, J. & Schäfer, K. 2001. Konstruieren im Stahlbetonbau. *Betonkalender 2001* Teil II: 311-492. Berlin: Ernst & Sohn.

Schlaich, J., Schäfer, K., Jennewein, M. 1987. Toward a Consistent Design of Structural Concrete. *PCI Journal* 32(3): 74-150.

Simo, J.C. & Hughes T.J.R. 1998. *Computational inelasticity*. New York: Springer.